# THE TRAP OF COMPLACENCY IN PREDICTING THE MAXIMUM

By J. du Toit and G. Peskir[1]

*University of the Witwatersrand and The University of Manchester*

Given a standard Brownian motion $B^\mu = (B_t^\mu)_{0 \le t \le T}$ with drift $\mu \in \mathbb{R}$ and letting $S_t^\mu = \max_{0 \le s \le t} B_s^\mu$ for $0 \le t \le T$, we consider the optimal prediction problem:

$$V = \inf_{0 \le \tau \le T} \mathsf{E}(B_\tau^\mu - S_T^\mu)^2$$

where the infimum is taken over all stopping times $\tau$ of $B^\mu$. Reducing the optimal prediction problem to a parabolic free-boundary problem we show that the following stopping time is optimal:

$$\tau_* = \inf\{t_* \le t \le T \mid b_1(t) \le S_t^\mu - B_t^\mu \le b_2(t)\}$$

where $t_* \in [0, T)$ and the functions $t \mapsto b_1(t)$ and $t \mapsto b_2(t)$ are continuous on $[t_*, T]$ with $b_1(T) = 0$ and $b_2(T) = 1/2\mu$. If $\mu > 0$, then $b_1$ is decreasing and $b_2$ is increasing on $[t_*, T]$ with $b_1(t_*) = b_2(t_*)$ when $t_* \neq 0$. Using local time-space calculus we derive a coupled system of nonlinear Volterra integral equations of the second kind and show that the pair of optimal boundaries $b_1$ and $b_2$ can be characterized as the unique solution to this system. This also leads to an explicit formula for $V$ in terms of $b_1$ and $b_2$. If $\mu \le 0$, then $t_* = 0$ and $b_2 \equiv +\infty$ so that $\tau_*$ is expressed in terms of $b_1$ only. In this case $b_1$ is decreasing on $[z_*, T]$ and increasing on $[0, z_*)$ for some $z_* \in [0, T)$ with $z_* = 0$ if $\mu = 0$, and the system of two Volterra equations reduces to one Volterra equation. If $\mu = 0$, then there is a closed form expression for $b_1$. This problem was solved in [*Theory Probab. Appl.* **45** (2001) 125–136] using the method of time change (i.e., change of variables). The

Received February 2005; revised March 2006.

[1]Research done at Network in Mathematical Physics and Stochastics (funded by the Danish National Research Foundation) and Centre for Analytical Finance (funded by the Danish Social Science Research Council).

method of time change cannot be extended to the case when $\mu \neq 0$ and the present paper settles the remaining cases using a different approach.

**1. Introduction.** Stopping a stochastic process as close as possible to its ultimate maximum is an undertaking of great practical and theoretical interest. Mathematical problems of this type may be referred to as *optimal prediction* problems. Variants of these problems have appeared in the past under different names (the optimal selection problem, the best choice problem, the secretary problem, the house selling problem) from where the older papers [1, 4, 7, 10] are interesting to consult. Most of this work has been done in the case of discrete time.

The case of continuous time has been studied in the recent papers [6] and [14] when the process is a standard Brownian motion. This hypothesis leads to an explicit solution of the problem using the method of time change. Motivated by wider applications, our aim in the present paper is to continue this study when the process is a standard Brownian motion with drift. It turns out that this extension is not only far from being routine, but also requires a different line of argument to be developed, which in turn is applicable to a broader class of diffusions and Markov processes.

The main objectives of the present paper may be described as follows. Given a standard Brownian motion $B^\mu = (B_t^\mu)_{0 \le t \le T}$ with drift $\mu \in \mathbb{R}$ and letting $S_t^\mu = \max_{0 \le s \le t} B_s^\mu$ for $0 \le t \le T$, we consider the optimal prediction problem:

$$(1.1) \qquad V = \inf_{0 \le \tau \le T} \mathsf{E}(B_\tau^\mu - S_T^\mu)^2$$

where the infimum is taken over all stopping times $\tau$ of $B^\mu$.

Reducing the optimal prediction problem (1.1) to a parabolic free-boundary problem we show that the following stopping time is optimal:

$$(1.2) \qquad \tau_* = \inf\{t_* \le t \le T \mid b_1(t) \le S_t^\mu - B_t^\mu \le b_2(t)\}$$

where $t_* \in [0, T)$ and the functions $t \mapsto b_1(t)$ and $t \mapsto b_2(t)$ are continuous on $[t_*, T]$ with $b_1(T) = 0$ and $b_2(T) = 1/2\mu$. If $\mu > 0$, then $b_1$ is decreasing and $b_2$ is increasing on $[t_*, T]$ with $b_1(t_*) = b_2(t_*)$ when $t_* \neq 0$. Using local time-space calculus (cf. [15, 16, 17]) we derive a coupled system of nonlinear Volterra integral equations of the second kind and show that the pair of optimal boundaries $b_1$ and $b_2$ can be characterized as the unique solution to this system. This also leads to an explicit formula for $V$ in terms of $b_1$ and $b_2$.

If $\mu \le 0$, then $t_* = 0$ and $b_2 \equiv +\infty$ so that $\tau_*$ is expressed in terms of $b_1$ only. In this case $b_1$ is decreasing on $[z_*, T]$ and increasing on $[0, z_*)$ for some $z_* \in [0, T)$ with $z_* = 0$ if $\mu = 0$, and the system of two Volterra equations reduces to one Volterra equation. If $\mu = 0$, then there is a closed form expression for $b_1$. This problem was solved in [6] using the method of time change (see also [14]). The method of time change cannot be extended to the case when $\mu \neq 0$ and the present paper settles the remaining cases using a different approach.

The continuation region of the problem turns out to be "humped" when $\mu < 0$. This is rather unexpected and indicates that the problem is strongly time dependent. The most surprising discovery revealed by the solution, however, is the existence of a nontrivial stopping region (a "black hole" as we call it) when $\mu > 0$. This fact is not only counter-intuitive but also has important practical implications. For example, in a growing economy where the appreciation rate of a stock price is strictly positive, any financial strategy based on optimal prediction of the ultimate maximum should be thoroughly re-examined in the light of this new phenomenon.

**2. The optimal prediction problem.** 1. Let $B = (B_t)_{t \geq 0}$ be a standard Brownian motion defined on a probability space $(\Omega, \mathcal{F}, \mathsf{P})$ where $B_0 = 0$ under $\mathsf{P}$. Set

$$B_t^\mu = B_t + \mu t \tag{2.1}$$

for $t \geq 0$ where $\mu \in \mathbb{R}$ is given and fixed. Then $B^\mu = (B_t^\mu)_{t \geq 0}$ is a standard Brownian motion with drift $\mu$. Define

$$S_t^\mu = \max_{0 \leq s \leq t} B_s^\mu \tag{2.2}$$

for $t \geq 0$. Then $S^\mu = (S_t^\mu)_{t \geq 0}$ is the maximum process associated with $B^\mu$.

2. Given $T > 0$ we consider the optimal prediction problem:

$$V = \inf_{0 \leq \tau \leq T} \mathsf{E}(B_\tau^\mu - S_T^\mu)^2 \tag{2.3}$$

where the infimum is taken over all stopping times $\tau$ of $B^\mu$ (the latter means that $\tau$ is a stopping time with respect to the natural filtration of $B^\mu$ that in turn is the same as the natural filtration of $B$ given by $\mathcal{F}_t^B = \sigma(B_s \mid 0 \leq s \leq t)$ for $t \in [0, T]$). The problem (2.3) consists of finding an optimal stopping time (at which the infimum is attained) and computing $V$ as explicitly as possible.

3. The identity (2.4) below reduces the optimal prediction problem (2.3) above [where the gain process $(B_t^\mu - S_T^\mu)_{0 \leq t \leq T}$ is not adapted to the natural filtration of $B^\mu$] to the optimal stopping problem (2.10) below (where the gain process is adapted). Similar arguments were exploited in [6] and [14] in the case when $\mu = 0$ in (2.3).

LEMMA 2.1. *The following identity holds:*

$$\mathsf{E}((S_T^\mu - B_t^\mu)^2 \mid \mathcal{F}_t^B) = (S_t^\mu - B_t^\mu)^2 + 2 \int_{S_t^\mu - B_t^\mu}^\infty z(1 - F^\mu(T - t, z)) \, dz \tag{2.4}$$

*for all $0 \leq t \leq T$ where*

$$
\begin{aligned}
F^\mu(T-t, z) &= \mathsf{P}(S^\mu_{T-t} \leq z) \\
&= \Phi\left(\frac{z - \mu(T-t)}{\sqrt{T-t}}\right) - e^{2\mu z}\Phi\left(\frac{-z - \mu(T-t)}{\sqrt{T-t}}\right)
\end{aligned}
\tag{2.5}
$$

*for $z \geq 0$.*

PROOF. By stationary independent increments of $B^\mu$ we have

$$
\begin{aligned}
\mathsf{E}((S^\mu_T - B^\mu_t)^2 | \mathcal{F}^B_t) \\
&= \mathsf{E}\left(\left(S^\mu_t + \left(\max_{t \leq s \leq T} B^\mu_s - S^\mu_t\right)^+ - B^\mu_t\right)^2 \bigg| \mathcal{F}^B_t\right) \\
&= \mathsf{E}\left(\left(S^\mu_t - B^\mu_t + \left(\max_{t \leq s \leq T} B^\mu_s - B^\mu_t - (S^\mu_t - B^\mu_t)\right)^+\right)^2 \bigg| \mathcal{F}^B_t\right) \\
&= (\mathsf{E}(x + (S^\mu_{T-t} - x)^+)^2)|_{x = S^\mu_t - B^\mu_t}
\end{aligned}
\tag{2.6}
$$

for $0 \leq t \leq T$ given and fixed. Integration by parts gives

$$
\begin{aligned}
\mathsf{E}(x + (S^\mu_{T-t} - x)^+)^2 &= \mathsf{E}(x^2 I(S^\mu_{T-t} \leq x)) + \mathsf{E}((S^\mu_{T-t})^2 I(S^\mu_{T-t} > x)) \\
&= x^2 \mathsf{P}(S^\mu_{T-t} \leq x) + \int_x^\infty z^2 F^\mu(T-t, dz) \\
&= x^2 F^\mu(T-t, x) + (z^2(F^\mu(T-t, z) - 1))|_x^\infty \\
&\quad + 2\int_x^\infty z(1 - F^\mu(T-t, z))\, dz \\
&= x^2 + 2\int_x^\infty z(1 - F^\mu(T-t, z))\, dz
\end{aligned}
\tag{2.7}
$$

for all $x \geq 0$. Combining (2.6) and (2.7) we get (2.4). The identity (2.5) is a well-known result of [3], page 397 and [13], page 526 completing the proof. □

4. Standard arguments based on the fact that each stopping time is the limit of a decreasing sequence of discrete stopping times imply that (2.4) extends as follows:

$$
\mathsf{E}((S^\mu_T - B^\mu_\tau)^2 | \mathcal{F}^B_\tau) = (S^\mu_\tau - B^\mu_\tau)^2 + 2\int_{S^\mu_\tau - B^\mu_\tau}^\infty z(1 - F^\mu(T-\tau, z))\, dz
\tag{2.8}
$$

for all stopping times $\tau$ of $B^\mu$ with values in $[0, T]$. Setting

$$
X_t = S^\mu_t - B^\mu_t
\tag{2.9}
$$

for $t \geq 0$ and taking expectations in (2.8) we find that the optimal prediction problem (2.3) is equivalent to the optimal stopping problem:

$$V = \inf_{0 \leq \tau \leq T} \mathsf{E}\left(X_\tau^2 + 2\int_{X_\tau}^\infty z(1 - F^\mu(T-\tau, z))\, dz\right) \tag{2.10}$$

where the infimum is taken over all stopping times $\tau$ of $X$ (upon recalling that the natural filtrations of $B^\mu$ and $X$ coincide). The process $X = (X_t)_{t \geq 0}$ is strong Markov so that (2.10) falls into the class of optimal stopping problems for Markov processes (cf. [19]). The structure of (2.10) is complicated since the gain process depends on time in a highly nonlinear way.

5. A successful treatment of (2.10) requires that the problem be extended so that the process $X$ can start at arbitrary points in the state space $[0, \infty)$. For this, recall that (cf. [5]) the following identity in law holds:

$$X \stackrel{\text{law}}{=} |Y| \tag{2.11}$$

where $|Y| = (|Y_t|)_{t \geq 0}$ and the process $Y = (Y_t)_{t \geq 0}$ is a unique strong solution to the stochastic differential equation:

$$dY_t = -\mu\, \text{sign}(Y_t)\, dt + dB_t \tag{2.12}$$

with $Y_0 = 0$. Moreover, it is known (cf. [5]) that under $Y_0 = x$ in (2.12) the process $|Y|$ has the same law as a Brownian motion with drift $-\mu$ started at $|x|$ and reflected at $0$. The infinitesimal operator of $|Y|$ acts on functions $f \in C_b^2([0, \infty))$ satisfying $f'(0+) = 0$ as $-\mu f'(x) + (1/2) f''(x)$. Since an optimal stopping time in (2.10) is the first entry time of the process to a closed set (this follows by general optimal stopping results and will be made more precise below) it is possible to replace the process $X$ in (2.10) by the process $|Y|$. On the other hand, since it is difficult to solve the equation (2.12) explicitly so that the dependence of $X$ on $x$ is clearly expressed, we will take a different route based on the following fact.

LEMMA 2.2. *The process $X^x = (X_t^x)_{t \geq 0}$ defined by*

$$X_t^x = x \vee S_t^\mu - B_t^\mu \tag{2.13}$$

*is Markov under $\mathsf{P}$ making $\mathsf{P}_x = \text{Law}(X^x \mid \mathsf{P})$ for $x \geq 0$ a family of probability measures on the canonical space $(C_+, \mathcal{B}(C_+))$ under which the coordinate process $X = (X_t)_{t \geq 0}$ is Markov with $\mathsf{P}_x(X_0 = x) = 1$.*

PROOF. Let $x \geq 0$, $t \geq 0$ and $h > 0$ be given and fixed. We then have

$$\begin{aligned}
X_{t+h}^x &= x \vee S_{t+h}^\mu - B_{t+h}^\mu \\
&= (x \vee S_t^\mu) \vee \left(\max_{t \leq s \leq t+h} B_s^\mu\right) - (B_{t+h}^\mu - B_t^\mu) - B_t^\mu \\
&= (x \vee S_t^\mu - B_t^\mu) \vee \left(\max_{t \leq s \leq t+h} B_s^\mu - B_t^\mu\right) - (B_{t+h}^\mu - B_t^\mu).
\end{aligned} \tag{2.14}$$

Hence by stationary independent increments of $B^\mu$ we get

(2.15) $\qquad \mathsf{E}(f(X^x_{t+h}) \mid \mathcal{F}^B_t) = \mathsf{E}(f(z \vee S^\mu_h - B^\mu_h))|_{z=X^x_t}$

for every bounded Borel function $f$. This shows that $X^x$ is a Markov process under $\mathsf{P}$. Moreover, the second claim follows from (2.15) by a basic transformation theorem for integrals upon using that the natural filtrations of $B$ and $X^x$ coincide. This completes the proof. $\square$

6. By means of Lemma 2.2 we can now extend the optimal stopping problem (2.10) where $X_0 = 0$ under $\mathsf{P}$ to the optimal stopping problem:

(2.16) $\quad V(t,x) = \inf_{0 \leq \tau \leq T-t} \mathsf{E}_{t,x}\left( X^2_{t+\tau} + 2\int_{X_{t+\tau}}^\infty z(1 - F^\mu(T-t-\tau, z))\,dz \right)$

where $X_t = x$ under $\mathsf{P}_{t,x}$ with $(t,x) \in [0,T] \times [0,\infty)$ given and fixed. The infimum in (2.16) is taken over all stopping times $\tau$ of $X$.

In view of the fact that $B^\mu$ has stationary independent increments, it is no restriction to assume that the process $X$ under $\mathsf{P}_{t,x}$ is explicitly given as

(2.17) $\qquad\qquad\qquad X^x_{t+s} = x \vee S^\mu_s - B^\mu_s$

under $\mathsf{P}$ for $s \in [0, T-t]$. Setting

(2.18) $\qquad\qquad\qquad R(t,z) = 1 - F^\mu(T-t, z)$

and introducing the gain function

(2.19) $\qquad\qquad G(t,x) = x^2 + 2\int_x^\infty z R(t,z)\,dz,$

we see that (2.16) can be written as follows:

(2.20) $\qquad\qquad V(t,x) = \inf_{0 \leq \tau \leq T-t} \mathsf{E}_{t,x}(G(t+\tau, X_{t+\tau}))$

for $(t,x) \in [0,T] \times [0,\infty)$.

7. The preceding analysis shows that the optimal prediction problem (2.3) reduces to solving the optimal stopping problem (2.20). Introducing the continuation set $C = \{(t,x) \in [0,T] \times [0,\infty) \mid V(t,x) < G(t,x)\}$ and the stopping set $D = \{(t,x) \in [0,T] \times [0,\infty) \mid V(t,x) = G(t,x)\}$, we may infer from general theory of optimal stopping for Markov processes (cf. [18]) that the optimal stopping time in (2.20) is given by

(2.21) $\qquad \tau_D = \inf\{0 \leq s \leq T-t \mid (t+s, X_{t+s}) \in D\}.$

It then follows using (2.9) that the optimal stopping time in (2.3) is given by

(2.22) $\qquad \tau_* = \inf\{0 \leq t \leq T \mid (t, S^\mu_t - B^\mu_t) \in D\}.$

The problems (2.20) and (2.3) are therefore reduced to determining $D$ and $V$ (outside $D$). We will see below that this task is complicated primarily because the gain function $G$ depends on time in a highly nonlinear way. The main aim of the paper is to present solutions to the problems formulated.

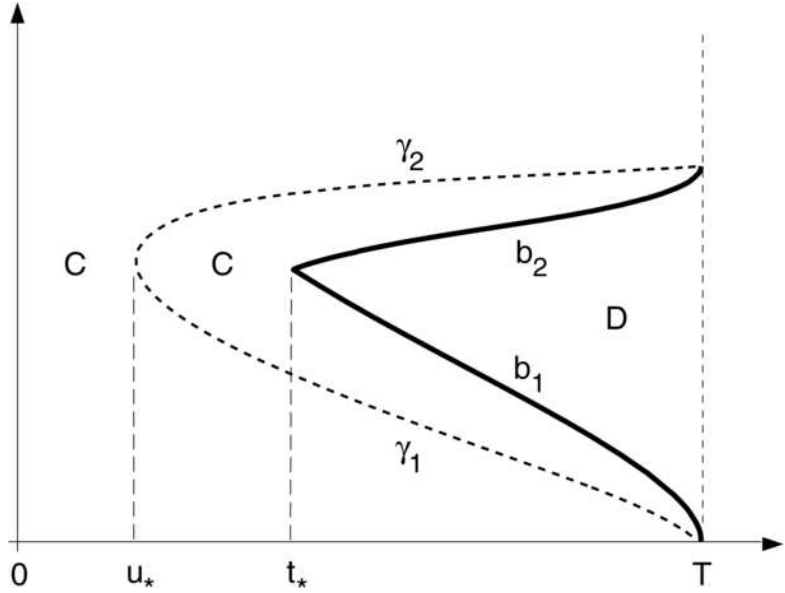

Fig. 1. (The "black-hole" effect.) *A computer drawing of the optimal stopping boundaries $b_1$ and $b_2$ when $\mu > 0$ is away from 0.*

**3. The free-boundary problem.** Consider the optimal stopping problem (2.20). Recall that the problem reduces to determining the stopping set $D$ and the value function $V$ outside $D$. It turns out that the shape of $D$ depends on the sign of $\mu$.

1. *The case $\mu > 0$*. It will be shown in the proof below that $D = \{(t,x) \in [t_*,T) \times [0,\infty) \mid b_1(t) \leq x \leq b_2(t)\} \cup \{(T,x) \mid x \in [0,\infty)\}$ where $t_* \in [0,T)$, the function $t \mapsto b_1(t)$ is continuous and decreasing on $[t_*,T]$ with $b_1(T) = 0$, and the function $t \mapsto b_2(t)$ is continuous and increasing on $[t_*,T]$ with $b_2(T) = 1/2\mu$. If $t_* \neq 0$, then $b_1(t_*) = b_2(t_*)$, and if $t_* = 0$, then $b_1(t_*) \leq b_2(t_*)$. We also have $b_1(t) < b_2(t)$ for all $t_* < t \leq T$. See Figures 1 and 2.

It follows that the optimal stopping time (2.21) can be written as follows:

$$(3.1) \qquad \tau_D = \inf\{t_* \leq t \leq T \mid b_1(t) \leq X_t \leq b_2(t)\}.$$

Inserting this expression into (2.20) and recalling that $C$ equals $D^c$ in $[0,T] \times [0,\infty)$, we can use Markovian arguments to formulate the following *free-boundary problem*:

$$(3.2) \quad V_t - \mu V_x + \tfrac{1}{2} V_{xx} = 0 \qquad \text{in } C,$$

$$(3.3) \qquad V(t,b_1(t)) = G(t,b_1(t)) \qquad \text{for } t_* \leq t \leq T,$$

$$(3.4) \qquad V(t,b_2(t)) = G(t,b_2(t)) \qquad \text{for } t_* \leq t \leq T,$$

$$(3.5) \qquad V_x(t,b_1(t)-) = G_x(t,b_1(t)) \qquad \text{for } t_* \leq t < T \text{ (smooth fit)},$$

(3.6) $\quad V_x(t, b_2(t)+) = G_x(t, b_2(t)) \quad$ for $t_* \leq t < T$ (smooth fit),

(3.7) $\quad V_x(t, 0+) = 0 \quad$ for $0 \leq t < T$ (normal reflection),

(3.8) $\quad V < G \quad$ in $C$,

(3.9) $\quad V = G \quad$ in $D$.

Note that the conditions (3.5)–(3.7) will be derived in the proof below while the remaining conditions are obvious.

2. *The case* $\mu \leq 0$. It will be seen in the proof below that $D = \{(t, x) \in [0, T) \times [0, \infty) \mid x \geq b_1(t)\} \cup \{(T, x) \mid x \in [0, \infty)\}$ where the continuous function $t \mapsto b_1(t)$ is decreasing on $[z_*, T]$ with $b_1(T) = 0$ and increasing on $[0, z_*)$ for some $z_* \in [0, T)$ (with $z_* = 0$ if $\mu = 0$). See Figures 3 and 4.

It follows that the optimal stopping time (2.21) can be written as follows:

(3.10) $\quad \tau_D = \inf\{0 \leq t \leq T \mid X_t \geq b_1(t)\}$.

Inserting this expression into (2.20) and recalling again that $C$ equals $D^c$ in $[0, T] \times [0, \infty)$, we can use Markovian arguments to formulate the following *free-boundary problem*:

(3.11) $\quad V_t - \mu V_x + \tfrac{1}{2} V_{xx} = 0 \quad$ in $C$,

(3.12) $\quad V(t, b_1(t)) = G(t, b_1(t)) \quad$ for $0 \leq t \leq T$,

(3.13) $\quad V_x(t, b_1(t)-) = G_x(t, b_1(t)) \quad$ for $0 \leq t < T$ (smooth fit),

(3.14) $\quad V_x(t, 0+) = 0 \quad$ for $0 \leq t < T$ (normal reflection),

(3.15) $\quad V < G \quad$ in $C$,

(3.16) $\quad V = G \quad$ in $D$.

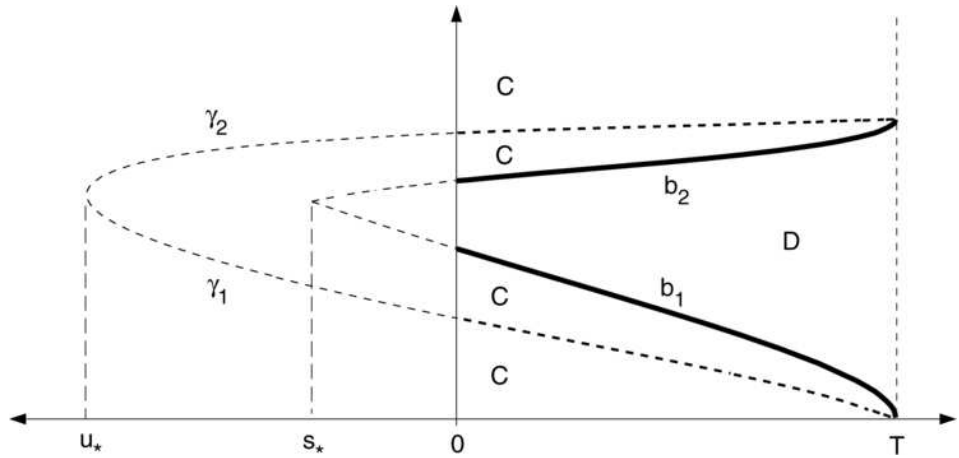

FIG. 2.  *A computer drawing of the optimal stopping boundaries $b_1$ and $b_2$ when $\mu \geq 0$ is close to* 0.

Note that the conditions (3.13) and (3.14) can be derived similarly to the conditions (3.5) and (3.7) above while the remaining conditions are obvious.

3. It will be clear from the proof below that the case $\mu \leq 0$ may be viewed as the case $\mu > 0$ with $b_2 \equiv \infty$ (and $t_* = 0$). This is in accordance with the facts that $b_2 \uparrow \infty$ as $\mu \downarrow 0$ and the point $s_* < T$ at which $b_1(s_*) = b_2(s_*)$ tends to $-\infty$ as $\mu \downarrow 0$. (Note that $t_*$ equals $s_* \vee 0$ and that extending the time interval $[0, T]$ to negative values in effect corresponds to enlarging the terminal value $T$ in the problem (2.20) above.) Since the case $\mu > 0$ is richer and more interesting we will only treat this case in complete detail. The case $\mu \leq 0$ can be dealt with analogously and most of the details will be omitted.

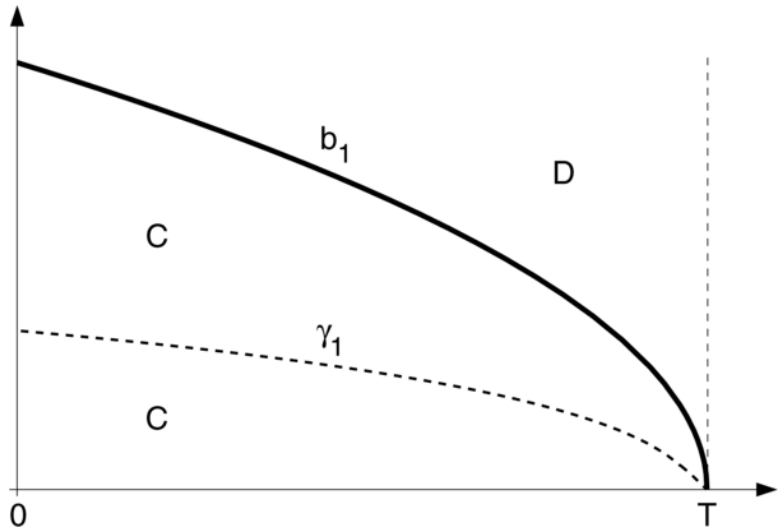

FIG. 3. *A computer drawing of the optimal stopping boundary $b_1$ when $\mu \leq 0$ is close to 0.*

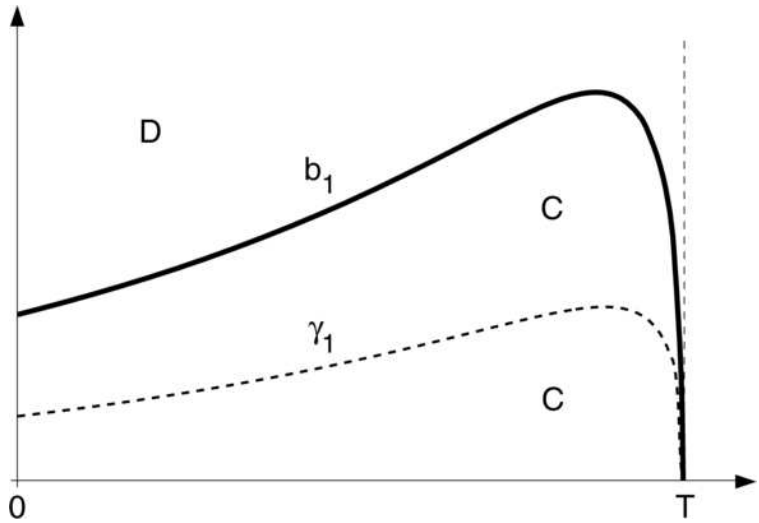

FIG. 4. *(The "hump" effect.) A computer drawing of the optimal stopping boundary $b_1$ when $\mu < 0$ is away from 0.*

4. It will follow from the result of Theorem 4.1 below that the free-boundary problem (3.2)–(3.9) characterizes the value function $V$ and the optimal stopping boundaries $b_1$ and $b_2$ in a unique manner. Motivated by wider application, however, our main aim will be to express $V$ in terms of $b_1$ and $b_2$ and show that $b_1$ and $b_2$ themselves satisfy a coupled system of nonlinear integral equations (which may then be solved numerically). Such an approach dates back to Kolodner [12] in a general context, while a particularly simple way of demonstrating these identities in the case of the American put problem (with one boundary) has been suggested in [2], [8] and [11]. The present problem, however, is in many ways different and substantially more complicated than the American put problem so that the proof will require novel arguments. We will nonetheless succeed in proving (as in [16] and [17] with one boundary) that the coupled system of nonlinear equations derived for $b_1$ and $b_2$ cannot have other solutions. The key argument in the proof relies upon a local time-space formula (see [15]). The analogous facts hold for the free-boundary problem (3.11)–(3.16) and the optimal stopping boundary $b_1$ (see Theorem 4.1 below).

**4. The result and proof.** 1. To solve the problems (2.3) and (2.20) let us introduce the function

$$(4.1) \qquad H = G_t - \mu G_x + \tfrac{1}{2} G_{xx}$$

on $[0, T] \times [0, \infty)$. A lengthy but straightforward calculation shows that

$$(4.2) \qquad \begin{aligned} H(t,x) &= (2\mu^2(T-t) - 2\mu x + 3)\Phi\!\left(\frac{x - \mu(T-t)}{\sqrt{T-t}}\right) \\ &\quad - 2\mu\sqrt{T-t}\,\varphi\!\left(\frac{x - \mu(T-t)}{\sqrt{T-t}}\right) \\ &\quad - e^{2\mu x}\Phi\!\left(\frac{-x - \mu(T-t)}{\sqrt{T-t}}\right) - 2(1 + \mu^2(T-t)) \end{aligned}$$

for $(t, x) \in [0, T] \times [0, \infty)$.

Let $P = \{(t,x) \in [0,T] \times [0,\infty) \mid H(t,x) \geq 0\}$ and $N = \{(t,x) \in [0,T] \times [0,\infty) \mid H(t,x) < 0\}$. A direct analysis based on (4.2) shows that in the case $\mu > 0$ we have $P = \{(t,x) \in [u_*, T] \times [0,\infty) \mid \gamma_1(t) \leq x \leq \gamma_2(t)\}$ where $u_* \in [0, T)$, the function $t \mapsto \gamma_1(t)$ is continuous and decreasing on $[u_*, T]$ with $\gamma_1(T) = 0$, and the function $t \mapsto \gamma_2(t)$ is continuous and increasing on $[u_*, T]$ with $\gamma_2(T) = 1/2\mu$. If $u_* \neq 0$, then $\gamma_1(u_*) = \gamma_2(u_*)$, and if $u_* = 0$, then $\gamma_1(u_*) \leq \gamma_2(u_*)$. We also have $\gamma_1(t) < \gamma_2(t)$ for all $u_* < t \leq T$. See Figures 1 and 2. Similarly, a direct analysis based on (4.2) shows that in the case $\mu \leq 0$ we have $P = \{(t, x) \in [0, T] \times [0, \infty) \mid x \geq \gamma_1(t)\}$ where the continuous function $t \mapsto \gamma_1(t)$ is decreasing on $[w_*, T]$ with $\gamma_1(T) = 0$ and increasing on $[0, w_*)$ for some $w_* \in [0, T)$ (with $w_* = 0$ if $\mu = 0$). See Figures 3 and 4.

2. Below we will make use of the following functions:

$$
\begin{aligned}
(4.3) \quad J(t, x) &= \mathsf{E}_x(G(T, X_{T-t})) \\
&= \int_0^\infty ds \int_{-\infty}^s db\, G(T, x \vee s - b) f(T-t, b, s),
\end{aligned}
$$

$$
\begin{aligned}
(4.4) \quad K(t, x, t+u, y, z) &= \mathsf{E}_x(H(t+u, X_u) I(y < X_u < z)) \\
&= \int_0^\infty ds \int_{-\infty}^s db\, H(t+u, x \vee s - b) \\
&\quad \times I(y < x \vee s - b < z) f(u, b, s),
\end{aligned}
$$

$$
\begin{aligned}
(4.5) \quad L(t, x, t+u, y) &= \mathsf{E}_x(H(t+u, X_u) I(X_u > y)) \\
&= \int_0^\infty ds \int_{-\infty}^s db\, H(t+u, x \vee s - b) \\
&\quad \times I(x \vee s - b > y) f(u, b, s),
\end{aligned}
$$

for $(t, x) \in [0, T] \times [0, \infty)$, $u \geq 0$ and $0 < y < z$, where $(b, s) \mapsto f(t, b, s)$ is the probability density function of $(B_t^\mu, S_t^\mu)$ under $\mathsf{P}$ given by

$$
(4.6) \quad f(t, b, s) = \sqrt{\frac{2}{\pi}} \frac{1}{t^{3/2}} (2s - b) e^{-(2s-b)^2/2t + \mu(b - \mu t/2)}
$$

for $t > 0$, $s \geq 0$ and $b \leq s$ (see, e.g., [9], page 368).

3. The main results of the paper may now be stated as follows.

THEOREM 4.1. *Consider the problems* (2.3) *and* (2.20). *We can then distinguish the following two cases:*

1. *The case* $\mu > 0$. *The optimal stopping boundaries in* (2.20) *can be characterized as the unique solution to the coupled system of nonlinear Volterra integral equations*

$$
\begin{aligned}
(4.7) \quad J(t, b_1(t)) &= G(t, b_1(t)) + \int_0^{T-t} K(t, b_1(t), t+u, b_1(t+u), b_2(t+u))\, du,
\end{aligned}
$$

$$
\begin{aligned}
(4.8) \quad J(t, b_2(t)) &= G(t, b_2(t)) + \int_0^{T-t} K(t, b_2(t), t+u, b_1(t+u), b_2(t+u))\, du,
\end{aligned}
$$

in the class of functions $t \mapsto b_1(t)$ and $t \mapsto b_2(t)$ on $[t_*, T]$ for $t_* \in [0, T)$ such that the function $t \mapsto b_1(t)$ is continuous and decreasing on $[t_*, T]$, the function $t \mapsto b_2(t)$ is continuous and increasing on $[t_*, T]$, and $\gamma_1(t) \leq b_1(t) < b_2(t) \leq \gamma_2(t)$ for all $t \in (t_*, T]$. The solutions $b_1$ and $b_2$ satisfy $b_1(T) = 0$ and $b_2(T) = 1/2\mu$, and the stopping time $\tau_D$ from (3.1) is optimal in (2.20). The stopping time (2.22) given by

$$(4.9) \qquad \tau_* = \inf\{0 \leq t \leq T \mid b_1(t) \leq S_t^\mu - B_t^\mu \leq b_2(t)\}$$

is optimal in (2.3). The value function $V$ from (2.20) admits the following representation:

$$(4.10) \quad V(t,x) = J(t,x) - \int_0^{T-t} K(t, x, t+u, b_1(t+u), b_2(t+u))\, du$$

for $(t, x) \in [0, T] \times [0, \infty)$. The value $V$ from (2.3) equals $V(0, 0)$ in (4.10).

2. The case $\mu \leq 0$. The optimal stopping boundary in (2.20) can be characterized as the unique solution to the nonlinear Volterra integral equation

$$(4.11) \quad J(t, b_1(t)) = G(t, b_1(t)) + \int_0^{T-t} L(t, b_1(t), t+u, b_1(t+u))\, du$$

in the class of continuous functions $t \mapsto b_1(t)$ on $[0, T]$ that are decreasing on $[z_*, T]$ and increasing on $[0, z_*)$ for some $z_* \in [0, T)$ and satisfy $b_1(t) \geq \gamma_1(t)$ for all $t \in [0, T]$. The solution $b_1$ satisfies $b_1(T) = 0$ and the stopping time $\tau_D$ from (3.10) is optimal in (2.20). The stopping time (2.22) given by

$$(4.12) \qquad \tau_* = \inf\{0 \leq t \leq T \mid S_t^\mu - B_t^\mu \geq b_1(t)\}$$

is optimal in (2.3). The value function $V$ from (2.20) admits the following representation:

$$(4.13) \qquad V(t,x) = J(t,x) - \int_0^{T-t} L(t, x, t+u, b_1(t+u))\, du$$

for $(t, x) \in [0, T] \times [0, \infty)$. The value $V$ from (2.3) equals $V(0, 0)$ in (4.13).

PROOF. The proof will be carried out in several steps. We will only treat the case $\mu > 0$ in complete detail. The case $\mu \leq 0$ can be dealt with analogously and details in this direction will be omitted. Thus we will assume throughout that $\mu > 0$ is given and fixed. We begin by invoking a result from general theory of optimal stopping for Markov processes.

1. We show that the stopping time $\tau_D$ in (2.21) is optimal in the problem (2.20). For this, recall that it is no restriction to assume that the process $X$ under $\mathsf{P}_{t,x}$ is given explicitly by (2.17) under $\mathsf{P}$. Since clearly $(t, x) \mapsto \mathsf{E}(G(t + \tau, X_\tau^x))$ is continuous (and thus upper semicontinuous) for each stopping time $\tau$, it follows that $(t, x) \mapsto V(t, x)$ is usc (recall that the infimum of usc functions defines a usc function). Since $(t, x) \mapsto G(t, x)$ is continuous (and thus lower semicontinuous) by general theory (cf. Corollary 2.7 in [18], Chapter 1, Section 2) it follows that $\tau_D$ is optimal in (2.20) as claimed. Note also that $C$ is open and $D$ is closed in $[0,T] \times [0,\infty)$.

2. The initial insight into the shape of $D$ is provided by stochastic calculus as follows. By Itô's formula we have

$$
\begin{aligned}
G(t+s, X_{t+s}) = G(t,x) &+ \int_0^s G_t(t+u, X_{t+u})\, du \\
&+ \int_0^s G_x(t+u, X_{t+u})\, dX_{t+u} \\
&+ \tfrac{1}{2} \int_0^s G_{xx}(t+u, X_{t+u})\, d\langle X, X\rangle_{t+u}
\end{aligned}
\tag{4.14}
$$

for $0 \le s \le T-t$ and $x \ge 0$ given and fixed. By the Itô–Tanaka formula [recalling (2.12) above] we have

$$
\begin{aligned}
X_t = |Y_t| &= x + \int_0^t \operatorname{sign}(Y_s) I(Y_s \ne 0)\, dY_s + \ell_t^0(Y) \\
&= x - \mu \int_0^t I(Y_s \ne 0)\, ds + \int_0^t \operatorname{sign}(Y_s) I(Y_s \ne 0)\, dB_s + \ell_t^0(Y)
\end{aligned}
\tag{4.15}
$$

where $\operatorname{sign}(0) = 0$ and $\ell_t^0(Y)$ is the local time of $Y$ at $0$ given by

$$
\ell_t^0(Y) = \mathsf{P} - \lim_{\varepsilon \downarrow 0} \frac{1}{2\varepsilon} \int_0^t I(-\varepsilon < Y_s < \varepsilon)\, ds
\tag{4.16}
$$

upon using that $d\langle Y, Y\rangle_s = ds$. It follows from (4.15) that

$$
dX_t = -\mu I(Y_t \ne 0)\, dt + \operatorname{sign}(Y_t) I(Y_t \ne 0)\, dB_t + d\ell_t^0(Y).
\tag{4.17}
$$

Inserting (4.17) into (4.14), using that $d\langle X, X\rangle_t = I(Y_t \ne 0)\, dt$ and $\mathsf{P}(Y_t = 0) = 0$, we get

$$
\begin{aligned}
G(t+s, X_{t+s}) = G(t,x) &+ \int_0^s (G_t - \mu G_x + \tfrac{1}{2} G_{xx})(t+u, X_{t+u})\, du \\
&+ \int_0^s G_x(t+u, X_{t+u}) \operatorname{sign}(Y_{t+u})\, dB_{t+u} \\
&+ \int_0^s G_x(t+u, X_{t+u})\, d\ell_{t+u}^0(Y) \\
&= G(t,x) + \int_0^s H(t+u, X_{t+u})\, du + M_s
\end{aligned}
\tag{4.18}
$$

where $H$ is given by (4.1) above and $M_s = \int_0^s G_x(t+u, X_{t+u}) \operatorname{sign}(Y_{t+u})\, dB_{t+u}$ is a continuous (local) martingale for $s \ge 0$. In the last identity in (4.18) we use that (quite remarkably) $G_x(t,0) = 0$ while $d\ell_{t+u}^0(Y)$ is concentrated at $0$ so that the final integral in (4.18) vanishes.

From the final expression in (4.18) we see that the initial insight into the shape of $D$ is gained by determining the sets $P$ and $N$ as introduced following (4.2) above. By considering the exit times from small balls in $[0, T) \times [0, \infty)$ and making use of (4.18) with the optional sampling theorem, we see that it is never optimal to stop in $N$. We thus conclude that $D \subseteq P$.

A deeper insight into the shape of $D$ is provided by the following arguments. Due to the fact that $P$ is bounded by $\gamma_1$ and $\gamma_2$ as described following (4.2) above, it is readily verified using (4.18) above and simple comparison arguments that for each $x \in (0, 1/2\mu)$ there exists $t = t(x) \in (0, T)$ close enough to $T$ such that every point $(x, u)$ belongs to $D$ for $u \in [t, T]$ (cf. [18], Chapter 8, Section 3). Note that this fact is fully in agreement with intuition since after starting at $(u, x)$ close to $(T, x)$ there will not be enough time to reach either of the favorable regions below $\gamma_1$ or above $\gamma_2$ to compensate for the loss incurred by strictly positive $H$ via (4.18). These arguments in particular show that $D \setminus \{(T, x) \mid x \in \mathbb{R}_+\}$ is nonempty.

The final insight into the shape of $D$ is obtained by the following fortunate fact:

(4.19) $\qquad\qquad\qquad t \mapsto H(t, x)$ is increasing on $[0, T]$

whenever $x \geq 0$. Indeed, this can be verified by a direct differentiation in (4.2) which yields

$$
\begin{aligned}
H_t(t, x) &= 2\left(\frac{x + \mu(T-t)}{(T-t)^{3/2}}\right)\varphi\left(\frac{x - \mu(T-t)}{\sqrt{T-t}}\right) \\
&\quad + 2\mu^2\left(1 - \Phi\left(\frac{x - \mu(T-t)}{\sqrt{T-t}}\right)\right)
\end{aligned}
$$
(4.20)

from where one sees that $H_t \geq 0$ on $[0, T) \times [0, \infty)$ upon recalling that $\mu > 0$ by assumption.

We next show that $(t_1, x) \in D$ implies that $(t_2, x) \in D$ whenever $0 \leq t_1 \leq t_2 \leq T$ and $x \geq 0$. For this, assume that $(t_2, x) \in C$ for some $t_2 \in (t_1, T)$. Let $\tau_* = \tau_D(t_2, x)$ denote the optimal stopping time for $V(t_2, x)$. Then by (4.18) and (4.19) using the optional sampling theorem we have

$$
\begin{aligned}
V(t_1, x) - G(t_1, x) &\leq \mathsf{E}(G(t_1 + \tau_*, X^x_{\tau_*})) - G(t_1, x) \\
&= \mathsf{E}\left(\int_0^{\tau_*} H(t_1 + u, X^x_u)\, du\right) \\
&\leq \mathsf{E}\left(\int_0^{\tau_*} H(t_2 + u, X^x_u)\, du\right) \\
&= \mathsf{E}(G(t_2 + \tau_*, X^x_{\tau_*}) - G(t_2, x)) \\
&= V(t_2, x) - G(t_2, x) < 0.
\end{aligned}
$$
(4.21)

Hence $(t_1, x)$ belongs to $C$, which is a contradiction. This proves the initial claim.

Finally we show that for $(t, x_1) \in D$ and $(t, x_2) \in D$ with $x_1 \le x_2$ in $(0, \infty)$ we have $(t, z) \in D$ for every $z \in [x_1, x_2]$. For this, fix $z \in (x_1, x_2)$ and let $\tau_* = \tau_D(t, z)$ denote the optimal stopping time for $V(t, z)$. Since $(u, x_1)$ and $(u, x_2)$ belong to $D$ for all $u \in [t, T]$ we see that $\tau_*$ must be smaller than or equal to the exit time from the rectangle $R$ with corners at $(t, x_1)$, $(t, x_2)$, $(T, x_1)$ and $(T, x_2)$. However, since $H > 0$ on $R$ we see from (4.18) upon using the optional sampling theorem that $V(t, z) > G(t, z)$. This shows that $(t, z)$ cannot belong to $C$, thus proving the initial claim.

Summarizing the facts derived above we can conclude that $D$ equals the set of all $(t, x)$ in $[t_*, T] \times [0, \infty)$ with $t_* \in [0, T)$ such that $b_1(t) \le x \le b_2(t)$, where the function $t \mapsto b_1(t)$ is decreasing on $[t_*, T]$ with $b_1(T) = 0$, the function $t \mapsto b_2(t)$ is increasing on $[t_*, T]$ with $b_2(T) = 1/2\mu$, and $\gamma_1(t) \le b_1(t) \le b_2(t) \le \gamma_2(t)$ for all $t \in [t_*, T]$. See Figures 1 and 2. It follows in particular that the stopping time $\tau_D$ from (3.1) is optimal in (2.20) and the stopping time from (4.9) is optimal in (2.3).

3. We show that $V$ is continuous on $[0, T] \times [0, \infty)$. For this, we will first show that $x \mapsto V(t, x)$ is continuous on $[0, \infty)$ uniformly over $t \in [0, T]$. Indeed, if $x < y$ in $[0, \infty)$ are given and fixed, we then have

$$
\begin{aligned}
&V(t, x) - V(t, y) \\
(4.22) \quad &= \inf_{0 \le \tau \le T-t} \mathsf{E}(G(t+\tau, X_\tau^x)) - \inf_{0 \le \tau \le T-t} \mathsf{E}(G(t+\tau, X_\tau^y)) \\
&\ge \inf_{0 \le \tau \le T-t} \mathsf{E}(G(t+\tau, X_\tau^x) - G(t+\tau, X_\tau^y))
\end{aligned}
$$

for all $t \in [0, T]$. It is easily verified that $x \mapsto G(t, x)$ is increasing so that $x \mapsto V(t, x)$ is increasing on $[0, \infty)$ for every $t \in [0, T]$. Hence it follows from (4.22) that

$$(4.23) \quad 0 \le V(t, y) - V(t, x) \le \sup_{0 \le \tau \le T-t} \mathsf{E}(G(t+\tau, X_\tau^y) - G(t+\tau, X_\tau^x))$$

for all $t \in [0, T]$. Using (2.19) we find

$$
\begin{aligned}
&G(t+\tau, X_\tau^y) - G(t+\tau, X_\tau^x) \\
&= (X_\tau^y)^2 - (X_\tau^x)^2 - 2\int_{X_\tau^x}^{X_\tau^y} zR(t+\tau, z)\,dt \\
(4.24) \quad &\le (X_\tau^y - X_\tau^x)(X_\tau^y + X_\tau^x + 2c) \\
&= (y \vee S_\tau^\mu - x \vee S_\tau^\mu)(y \vee S_\tau^\mu - B_\tau^\mu + x \vee S_\tau^\mu - B_\tau^\mu + 2c) \\
&\le (y - x)Z
\end{aligned}
$$

where $c = \sup_{z \geq 0} zR(t,z) \leq \mathsf{E}(S^\mu_{T-t}) \leq \mathsf{E}(S^\mu_T) < \infty$ by Markov's inequality and $Z = 2(y+1) + 4\max_{0 \leq t \leq T} |B^\mu_t| + 2c$ belongs to $L^1(\mathsf{P})$. From (4.23) and (4.24) we find

$$(4.25) \qquad 0 \leq V(t,y) - V(t,x) \leq (y-x)\mathsf{E}(Z)$$

for all $t \in [0,T]$ implying that $x \mapsto V(t,x)$ is continuous on $[0,\infty)$ uniformly over $t \in [0,T]$.

To complete the proof of the initial claim it is sufficient to show that $t \mapsto V(t,x)$ is continuous on $[0,T]$ for each $x \in [0,\infty)$ given and fixed. For this, fix $x$ in $[0,\infty)$ and $t_1 < t_2$ in $[0,T]$. Let $\tau_1 = \tau_D(t_1,x)$ and $\tau_2 = \tau_D(t_2,x)$ be optimal for $V(t_1,x)$ and $V(t_2,x)$, respectively. Setting $\tau_1^\varepsilon = \tau_1 \wedge (T-t_2)$ with $\varepsilon = t_2 - t_1$ we have

$$(4.26) \quad \begin{aligned} &\mathsf{E}(G(t_2+\tau_2, X^x_{\tau_2}) - G(t_1+\tau_2, X^x_{\tau_2})) \\ &\qquad \leq V(t_2,x) - V(t_1,x) \\ &\qquad \leq \mathsf{E}(G(t_2+\tau_1^\varepsilon, X^x_{\tau_1^\varepsilon}) - G(t_1+\tau_1, X^x_{\tau_1})). \end{aligned}$$

Note that we have

$$(4.27) \qquad G_t(t,x) = -2\int_x^\infty z f^\mu_{T-t}(z)\,dz$$

where $f^\mu_{T-t}(z) = (dF^\mu_{T-t}/dz)(z)$ so that

$$(4.28) \qquad |G_t(t,x)| \leq 2\int_0^\infty z f^\mu_{T-t}(z)\,dz = 2\mathsf{E}(S^\mu_{T-t}) \leq 2\mathsf{E}(S^\mu_T)$$

for all $t \in [0,T]$. Hence setting $\beta = 2\mathsf{E}(S^\mu_T)$ by the mean value theorem we get

$$(4.29) \qquad |G(u_2,x) - G(u_1,x)| \leq \beta(u_2 - u_1)$$

for all $u_1 < u_2$ in $[0,T]$. Using (4.29) in (4.26) upon subtracting and adding $G(t_1+\tau_1, X^x_{\tau_1^\varepsilon})$ we obtain

$$(4.30) \quad \begin{aligned} -\beta(t_2-t_1) &\leq V(t_2,x) - V(t_1,x) \\ &\leq 2\beta(t_2-t_1) + \mathsf{E}(G(t_1+\tau_1, X^x_{\tau_1^\varepsilon}) - G(t_1+\tau_1, X^x_{\tau_1})). \end{aligned}$$

Note that we have

$$(4.31) \qquad G_x(t,x) = 2xF^\mu_{T-t}(x) \leq 2x$$

so that the mean value theorem implies

$$(4.32) \quad \begin{aligned} |G(t_1+\tau_1, X^x_{\tau_1^\varepsilon}) - G(t_1+\tau_1, X^x_{\tau_1})| &= |G_x(t_1+\tau_1, \xi)||X^x_{\tau_1^\varepsilon} - X^x_{\tau_1}| \\ &\leq 2(X^x_{\tau_1^\varepsilon} \vee X^x_{\tau_1})|X^x_{\tau_1^\varepsilon} - X^x_{\tau_1}| \end{aligned}$$

where $\xi$ lies between $X^x_{\tau^\varepsilon_1}$ and $X^x_{\tau_1}$. Since $X^x_\tau$ is dominated by $x + 2\max_{0 \le t \le T} |B^\mu_t|$ which belongs to $L^1(\mathsf{P})$ for every stopping time $\tau$, letting $t_2 - t_1 \to 0$ and using that $\tau^\varepsilon_1 - \tau_1 \to 0$ we see from (4.30) and (4.32) that $V(t_2, x) - V(t_1, x) \to 0$ by dominated convergence. This shows that $t \mapsto V(t, x)$ is continuous on $[0, T]$ for each $x \in [0, \infty)$, and thus $V$ is continuous on $[0, T] \times [0, \infty)$ as claimed. Standard arguments based on the strong Markov property and classic results from PDEs show that $V$ is $C^{1,2}$ on $C$ and satisfies (3.2). These facts will be freely used below.

4. We show that $x \mapsto V(t, x)$ is differentiable at $b_i(t)$ for $i = 1, 2$ and that $V_x(t, b_i(t)) = G_x(t, b_i(t))$ for $t \in [t_*, T)$. For this, fix $t \in [t_*, T)$ and set $x = b_2(t)$ [the case $x = b_1(t)$ can be treated analogously]. We then have

$$(4.33) \qquad \frac{V(t, x+\varepsilon) - V(t, x)}{\varepsilon} \le \frac{G(t, x+\varepsilon) - G(t, x)}{\varepsilon}$$

for all $\varepsilon > 0$. Letting $\varepsilon \downarrow 0$ in (4.33) we find

$$(4.34) \qquad \limsup_{\varepsilon \downarrow 0} \frac{V(t, x+\varepsilon) - V(t, x)}{\varepsilon} \le G_x(t, x).$$

Let $\tau_\varepsilon = \tau_D(t, x + \varepsilon)$ be optimal for $V(t, x + \varepsilon)$. Then by the mean value theorem we have

$$(4.35) \qquad \begin{aligned} \frac{V(t, x+\varepsilon) - V(t, x)}{\varepsilon} &\ge \frac{1}{\varepsilon}(\mathsf{E}(G(t + \tau_\varepsilon, X^{x+\varepsilon}_{\tau_\varepsilon})) - \mathsf{E}(G(t + \tau_\varepsilon, X^x_{\tau_\varepsilon}))) \\ &= \frac{1}{\varepsilon}\mathsf{E}(G_x(t + \tau_\varepsilon, \xi_\varepsilon)(X^{x+\varepsilon}_{\tau_\varepsilon} - X^x_{\tau_\varepsilon})) \end{aligned}$$

where $\xi_\varepsilon$ lies between $X^x_{\tau_\varepsilon}$ and $X^{x+\varepsilon}_{\tau_\varepsilon}$. Using that $t \mapsto b_2(t)$ is increasing and that $t \mapsto \lambda t$ is a lower function for $B$ at $0+$ for every $\lambda \in \mathbb{R}$, it is possible to verify that $\tau_\varepsilon \to 0$ as $\varepsilon \downarrow 0$. Hence it follows that $\xi_\varepsilon \to x$ as $\varepsilon \downarrow 0$ so that $G_x(t + \tau_\varepsilon, \xi_\varepsilon) \to G_x(t, x)$ as $\varepsilon \downarrow 0$. Moreover, using (4.31) we find

$$(4.36) \qquad \begin{aligned} G_x(t + \tau_\varepsilon, \xi_\varepsilon) &\le 2\xi_\varepsilon \le 2 X^{x+\varepsilon}_{\tau_\varepsilon} = 2((x+\varepsilon) \vee S^\mu_{\tau_\varepsilon} - B^\mu_{\tau_\varepsilon}) \\ &\le 2\Big(x + \varepsilon + 2 \max_{0 \le t \le T} |B^\mu_t|\Big) \end{aligned}$$

where the final expression belongs to $L^1(\mathsf{P})$ (recall also that $G_x \ge 0$). Finally, we have

$$(4.37) \qquad \frac{1}{\varepsilon}(X^{x+\varepsilon}_{\tau_\varepsilon} - X^x_{\tau_\varepsilon}) = \frac{1}{\varepsilon}((x+\varepsilon) \vee S^\mu_{\tau_\varepsilon} - x \vee S^\mu_{\tau_\varepsilon}) \to 1$$

when $\varepsilon \downarrow 0$ as well as

$$(4.38) \qquad 0 \le \frac{1}{\varepsilon}(X^{x+\varepsilon}_{\tau_\varepsilon} - X^x_{\tau_\varepsilon}) \le 1$$

for all $\varepsilon > 0$. Letting $\varepsilon \downarrow 0$ in (4.35) and using (4.36)–(4.38), we may conclude that

$$\liminf_{\varepsilon \downarrow 0} \frac{V(t, x+\varepsilon) - V(t,x)}{\varepsilon} \geq G_x(t,x) \tag{4.39}$$

by dominated convergence. Combining (4.34) and (4.39) we see that $x \mapsto V(t,x)$ is differentiable at $b_2(t)$ with $V_x(t, b_2(t)) = G_x(t, b_2(t))$ as claimed. Analogously one finds that $x \mapsto V(t,x)$ is differentiable at $b_1(t)$ with $V_x(t, b_1(t)) = G_x(t, b_1(t))$ and further details of this derivation will be omitted.

A small modification of the proof above shows that $x \mapsto V(t,x)$ is $C^1$ at $b_2(t)$. Indeed, let $\tau_\delta = \tau_D(t, x+\delta)$ be optimal for $V(t, x+\delta)$ where $\delta > 0$ is given and fixed. Instead of (4.33) above we have by the mean value theorem that

$$\begin{aligned}
\frac{V(t, x+\delta+\varepsilon) - V(t, x+\delta)}{\varepsilon} \\
\leq \frac{1}{\varepsilon}(\mathsf{E}(G(t+\tau_\delta, X_{\tau_\delta}^{x+\delta+\varepsilon})) - \mathsf{E}(G(t+\tau_\delta, X_{\tau_\delta}^{x+\delta})))) \\
= \frac{1}{\varepsilon}\mathsf{E}(G_x(t+\tau_\delta, \eta_\varepsilon)(X_{\tau_\delta}^{x+\delta+\varepsilon} - X_{\tau_\delta}^{x+\delta}))
\end{aligned} \tag{4.40}$$

where $\eta_\varepsilon$ lies between $X_{\tau_\delta}^{x+\delta}$ and $X_{\tau_\delta}^{x+\delta+\varepsilon}$ for $\varepsilon > 0$. Clearly $\eta_\varepsilon \to X_{\tau_\delta}^{x+\delta}$ as $\varepsilon \downarrow 0$. Letting $\varepsilon \downarrow 0$ in (4.40) and using the same arguments as in (4.36)–(4.38) we can conclude that

$$V_x(t, x+\delta) \leq \mathsf{E}(G_x(t+\tau_\delta, X_{\tau_\delta}^{x+\delta})). \tag{4.41}$$

Moreover, in exactly the same way as in (4.35)–(4.39) we find that the reverse inequality in (4.41) also holds, so that we have

$$V_x(t, x+\delta) = \mathsf{E}(G_x(t+\tau_\delta, X_{\tau_\delta}^{x+\delta})). \tag{4.42}$$

Letting $\delta \downarrow 0$ in (4.42), recalling that $\tau_\delta \to 0$ and using the same arguments as in (4.36), we find by dominated convergence that

$$\lim_{\delta \downarrow 0} V_x(t, x+\delta) = G_x(t,x) = V_x(t,x). \tag{4.43}$$

Thus $x \mapsto V(t,x)$ is $C^1$ at $b_2(t)$ as claimed. Similarly one finds that $x \mapsto V(t,x)$ is $C^1$ at $b_1(t)$ with $V_x(t, b_1(t)+) = G_x(t, b_1(t))$ and further details of this derivation will be omitted. This establishes the smooth fit conditions (3.5), (3.6) and (3.13) above.

5. We show that $t \mapsto b_1(t)$ and $t \mapsto b_2(t)$ are continuous on $[t_*, T]$. Again we only consider the case of $b_2$ in detail, since the case of $b_1$ can be treated similarly. Note that the same proof also shows that $b_2(T-) = 1/2\mu$ and that $b_1(T-) = 0$.

Let us first show that $b_2$ is right-continuous. For this, fix $t \in [t_*, T)$ and consider a sequence $t_n \downarrow t$ as $n \to \infty$. Since $b_2$ is increasing, the right-hand limit $b_2(t+)$ exists. Because $(t_n, b_2(t_n))$ belongs to $D$ for all $n \geq 1$, and $D$ is closed, it follows that $(t, b_2(t+))$ belongs to $D$. Hence by (3.1) we may conclude that $b_2(t+) \leq b_2(t)$. Since the fact that $b_2$ is increasing gives the reverse inequality, it follows that $b_2$ is right-continuous as claimed.

Let us next show that $b_2$ is left-continuous. For this, suppose that there exists $t \in (t_*, T)$ such that $b_2(t-) < b_2(t)$. Fix a point $x \in (b_2(t-), b_2(t))$ and note by (3.6) that we have

$$
(4.44) \quad \begin{aligned} V(s,x) - G(s,x) \\ = \int_{b_2(s)}^{x} \int_{b_2(s)}^{y} (V_{xx}(s,z) - G_{xx}(s,z)) \, dz \, dy \end{aligned}
$$

for any $s \in (t_*, t)$. By (3.2) and (4.1) we find that

$$
(4.45) \quad \tfrac{1}{2}(V_{xx} - G_{xx}) = G_t - V_t + \mu(V_x - G_x) - H.
$$

From (4.21) we derive the key inequality

$$
(4.46) \quad V_t(t,x) \geq G_t(t,x)
$$

for all $(t,x) \in [0,T) \times [0,\infty)$. Inserting (4.45) into (4.44) and using (4.46) and (3.4) we find

$$
(4.47) \quad \begin{aligned} V(s,x) - G(s,x) &\leq \int_{b_2(s)}^{x} \int_{b_2(s)}^{y} 2(\mu(V_x - G_x)(s,z) - H(s,z)) \, dz \, dy \\ &= \int_{b_2(s)}^{x} 2\mu(V(s,y) - G(s,y)) \, dy \\ &\quad - \int_{b_2(s)}^{x} \int_{b_2(s)}^{y} 2H(s,z) \, dz \, dy \\ &\leq - \int_{b_2(s)}^{x} \int_{b_2(s)}^{y} 2H(s,z) \, dz \, dy \end{aligned}
$$

for any $s \in (t_*, t)$. From the properties of the function $\gamma_2$ it follows that there exists $s_* < t$ close enough to $t$ such that $(s, z)$ belongs to $P$ for all $s \in [s_*, t)$ and $z \in [b_2(s), x]$. Moreover, since $H$ is continuous and thus attains its infimum on a compact set, it follows that $2H(s,z) \geq m > 0$ for all $s \in [s_*, t)$ and $z \in [b_2(s), x]$. Using this fact in (4.47) we get

$$
(4.48) \quad V(s,x) - G(s,x) \leq -m \frac{(x - b_2(s))^2}{2} < 0
$$

for all $s \in [s_*, t)$. Letting $s \uparrow t$ in (4.48) we conclude that $V(t, x) < G(t, x)$ violating the fact that $(t, x) \in D$. This shows that $b_2$ is left-continuous and thus continuous. The continuity of $b_1$ is proved analogously.

6. We show that the normal reflection condition (3.7) holds. For this, note first since $x \mapsto V(t, x)$ is increasing on $[0, \infty)$ that $V_x(t, 0+) \geq 0$ for all $t \in [0, T)$ (note that the limit exists since $V$ is $C^{1,2}$ on $C$). Suppose that there exists $t \in [0, T)$ such that $V_x(t, 0+) > 0$. Recalling that $V$ is $C^{1,2}$ on $C$ so that $t \mapsto V_x(t, 0+)$ is continuous on $[0, T)$, we see that there exists $\delta > 0$ such that $V_x(s, 0+) \geq \varepsilon > 0$ for all $s \in [t, t + \delta]$ with $t + \delta < T$. Setting $\tau_\delta = \tau_D \wedge \delta$ it follows by the Itô–Tanaka formula [as in (4.18) above] upon using (3.2) and the optional sampling theorem [recall (4.41) and (4.31) for the latter] that we have

$$\begin{aligned}
\mathsf{E}_{t,0}&(V(t + \tau_\delta, X_{t+\tau_\delta})) \\
(4.49) \qquad &= V(t, 0) + \mathsf{E}_{t,0}\bigg(\int_0^{\tau_\delta} V_x(t + u, X_{t+u}) \, d\ell_{t+u}^0(Y)\bigg) \\
&\geq V(t, 0) + \varepsilon \mathsf{E}_{t,0}(\ell_{t+\tau_\delta}^0(Y)).
\end{aligned}$$

Since $(V(t + s \wedge \tau_D, X_{t+s \wedge \tau_D}))_{0 \leq s \leq T-t}$ is a martingale under $\mathsf{P}_{t,0}$ by general theory of optimal stopping for Markov processes (see, e.g., [18]) we see from (4.49) that $\mathsf{E}_{t,0}(\ell_{t+\tau_\delta}^0(Y))$ must be equal to 0. Since, however, properties of the local time clearly exclude this, we must have $V(t, 0+)$ equal to 0 as claimed in (3.7) above.

7. We show that $V$ is given by the formula (4.10) and that $b_1$ and $b_2$ solve the system (4.7)–(4.8). For this, note that by (3.2) and (4.46) we have $\frac{1}{2} V_{xx} = -V_t + \mu V_x \leq -G_t + \mu V_x$ in $C$. It is easily verified using (4.31) and (4.41) that $V_x(t, x) \leq M/2\mu$ for all $t \in [0, T)$ and all $x \in [0, (1/2\mu) + 1]$ with some $M > 0$ large enough. Using this inequality in the previous inequality we get $V_{xx} \leq -G_t + M$ in $A = C \cap ([0, T) \times [0, (1/2\mu) + 1])$. Setting $h(t, x) = \int_0^x \int_0^y (-G_t(t, z) + M) \, dz \, dy$ we easily see that $h$ is $C^{1,2}$ on $[0, T) \times [0, \infty)$ and that $h_{xx} = -G_t + M$. Thus the previous inequality reads $V_{xx} \leq h_{xx}$ in $A$, and setting $F = V - h$ we see that $x \mapsto F(t, x)$ is concave on $[0, b_1(t)]$ and $[b_2(t), (1/2\mu) + 1]$ for $t \in [t_*, T)$. We also see that $F$ is $C^{1,2}$ on $C$ and $D^o = \{(t, x) \in [t_*, T) \times [0, \infty) \mid b_1(t) < x < b_2(t)\}$ since both $V$ and $G$ are so. Moreover, it is also clear that $F_t - \mu F_x + \frac{1}{2} F_{xx}$ is locally bounded on $C \cup D^o$ in the sense that the function is bounded on $K \cap (C \cup D^o)$ for each compact set $K$ in $[0, T) \times [0, \infty)$. Finally, we also see using (3.5) and (3.6) that $t \mapsto F_x(t, b_i(t) \mp) = V_x(t, b_i(t) \mp) - h_x(t, b_i(t) \mp) = G_x(t, b_i(t)) - h_x(t, b_i(t))$ is continuous on $[t_*, T)$ since $b_i$ is continuous for $i = 1, 2$.

Since the previous conditions are satisfied we know that the local time-space formula (cf. Theorem 3.1 and Remark 2.3 in [15]) can be applied to $F(t + s, X_{t+s})$. Since $h$ is $C^{1,2}$ on $[0, T) \times [0, \infty)$ we know that the Itô–Tanaka formula can be applied to $h(t + s, X_{t+s})$ as in (4.18) above [upon noting that $h_x(t,0+) = 0$]. Adding the two formulae, using in the former that $F_x(t,0+) = -h_x(t,0+) = 0$ since $V_x(t,0+) = 0$ by (3.7) above, we get

$$
\begin{aligned}
V(t+s, X_{t+s}) &= V(t,x) \\
&+ \int_0^s (V_t - \mu V_x + \tfrac{1}{2} V_{xx})(t+u, X_{t+u}) \\
&\quad \times I(X_{t+u} \notin \{b_1(t+u), b_2(t+u)\})\, du \\
&+ \int_0^s V_x(t+u, X_{t+u}) \operatorname{sign}(Y_{t+u}) \\
&\quad \times I(X_{t+u} \notin \{b_1(t+u), b_2(t+u)\})\, dB_{t+u} \\
&+ \sum_{i=1}^2 \int_0^s (V_x(t+u, X_{t+u}+) - V_x(t+u, X_{t+u}-)) \\
&\quad \times I(X_{t+u} = b_i(t+u))\, d\ell_{t+u}^{b_i}(X)
\end{aligned}
$$
(4.50)

for $t \in [0,T)$ and $x \in [0,\infty)$. Making use of (3.2) and (3.9) in the first integral and (3.5) and (3.6) in the final integral (which consequently vanishes), we obtain

$$
\begin{aligned}
V(t+s, X_{t+s}) &= V(t,x) \\
&+ \int_0^s H(t+u, X_{t+u}) \\
&\quad \times I(b_1(t+u) < X_{t+u} < b_2(t+u))\, du + M_s
\end{aligned}
$$
(4.51)

for $t \in [0,T)$ and $x \in [0,\infty)$ where $M_s = \int_0^s V_x(t+u, X_{t+u})\, dB_{t+u}$ is a continuous (local) martingale for $s \ge 0$.

Setting $s = T-t$, using that $V(T,x) = G(T,x)$ for all $x \ge 0$ and taking the $\mathsf{P}_{t,x}$-expectation in (4.51), we find by the optional sampling theorem that

$$
\begin{aligned}
V(t,x) &= \mathsf{E}_{t,x}(G(T, X_T)) \\
&\quad - \int_0^{T-t} \mathsf{E}_{t,x}(H(t+u, X_{t+u}) \\
&\quad \times I(b_1(t+u) < X_{t+u} < b_2(t+u)))\, du
\end{aligned}
$$
(4.52)

for $t \in [0,T)$ and $x \in [0,\infty)$. Making use of (4.3) and (4.4) we see that (4.52) is the formula (4.10). Moreover, inserting $x = b_i(t)$ in (4.52) and using that $V(t, b_i(t)) = G(t, b_i(t))$ for $i = 1, 2$ we see that $b_1$ and $b_2$ satisfy the system (4.7)–(4.8) as claimed.

8. We show that $b_1$ and $b_2$ are the unique solution to the system (4.7)–(4.8) in the class of continuous functions $t \mapsto b_1(t)$ and $t \mapsto b_2(t)$ on $[t_*, T]$ for $t_* \in [0,T)$ such that $\gamma_1(t) \le b_1(t) < b_2(t) \le \gamma_2(t)$ for all $t \in (t_*, T]$. Note that there is no need to assume that $b_1$ is decreasing and $b_2$ is increasing as established above. The proof of uniqueness will be presented in the final three steps of the main proof below.

9. Let $c_1 : [t_*, T] \to \mathbb{R}$ and $c_2 : [t_*, T] \to \mathbb{R}$ be a solution to the system (4.7)–(4.8) for $t_* \in [0, T)$ such that $c_1$ and $c_2$ are continuous and satisfy $\gamma_1(t) \leq c_1(t) < c_2(t) \leq \gamma_2(t)$ for all $t \in (t_*, T]$. We need to show that these $c_1$ and $c_2$ must then be equal to the optimal stopping boundaries $b_1$ and $b_2$, respectively.

Motivated by the derivation (4.50)–(4.52) which leads to the formula (4.10), let us consider the function $U^c : [0, T) \times [0, \infty) \to \mathbb{R}$ defined as follows:

$$
\begin{aligned}
U^c(t, x) = {}& \mathsf{E}_{t,x}(G(T, X_T)) \\
& - \int_0^{T-t} \mathsf{E}_{t,x}(H(t+u, X_{t+u}) \\
& \qquad \times I(c_1(t+u) < X_{t+u} < c_2(t+u))) \, du
\end{aligned}
\tag{4.53}
$$

for $(t, x) \in [0, T) \times [0, \infty)$. In terms of (4.3) and (4.4) note that $U^c$ is explicitly given by

$$
U^c(t, x) = J(t, x) - \int_0^{T-t} K(t, x, t+u, c_1(t+u), c_2(t+u)) \, du \tag{4.54}
$$

for $(t, x) \in [0, T) \times [0, \infty)$. Observe that the fact that $c_1$ and $c_2$ solve the system (4.7)–(4.8) means exactly that $U^c(t, c_i(t)) = G(t, c_i(t))$ for all $t \in [t_*, T]$ and $i = 1, 2$. We will moreover show that $U^c(t, x) = G(t, x)$ for all $x \in [c_1(t), c_2(t)]$ with $t \in [t_*, T]$. This is the key point in the proof (cf. [16] and [17]) that can be derived using martingale arguments as follows.

If $X = (X_t)_{t \geq 0}$ is a Markov process (with values in a general state space) and we set $F(t, x) = \mathsf{E}_x(G(X_{T-t}))$ for a (bounded) measurable function $G$ with $\mathsf{P}(X_0 = x) = 1$, then the Markov property of $X$ implies that $F(t, X_t)$ is a martingale under $\mathsf{P}_x$ for $0 \leq t \leq T$. Similarly, if we set $F(t, x) = \mathsf{E}_x(\int_0^{T-t} H(X_s) \, ds)$ for a (bounded) measurable function $H$ with $\mathsf{P}(X_0 = x) = 1$, then the Markov property of $X$ implies that $F(t, X_t) + \int_0^t H(X_s) \, ds$ is a martingale under $\mathsf{P}_x$ for $0 \leq t \leq T$.

Combining the two martingale facts applied to the time-space Markov process $(t+s, X_{t+s})$ instead of $X_s$, we find that

$$
\begin{aligned}
& U^c(t+s, X_{t+s}) \\
& - \int_0^s H(t+u, X_{t+u}) I(c_1(t+u) < X_{t+u} < c_2(t+u)) \, du
\end{aligned}
\tag{4.55}
$$

is a martingale under $\mathsf{P}_{t,x}$ for $0 \leq s \leq T - t$. We may thus write

$$
\begin{aligned}
& U^c(t+s, X_{t+s}) \\
& - \int_0^s H(t+u, X_{t+u}) I(c_1(t+u) < X_{t+u} < c_2(t+u)) \, du \\
& = U^c(t, x) + N_s
\end{aligned}
\tag{4.56}
$$

where $(N_s)_{0 \leq s \leq T-t}$ is a martingale under $\mathsf{P}_{t,x}$. On the other hand, we know from (4.18) that

(4.57) $\quad G(t+s, X_{t+s}) = G(t,x) + \int_0^s H(t+u, X_{t+u})\,du + M_s$

where $M_s = \int_0^s G_x(t+u, X_{t+u})\,\text{sign}(Y_{t+u})\,dB_{t+u}$ is a continuous (local) martingale under $\mathsf{P}_{t,x}$ for $0 \leq s \leq T-t$.

For $x \in [c_1(t), c_2(t)]$ with $t \in [t_*, T]$ given and fixed, consider the stopping time

(4.58) $\quad \sigma_c = \inf\{0 \leq s \leq T-t \mid X_{t+s} \leq c_1(t+s) \text{ or } X_{t+s} \geq c_2(t+s)\}$

under $\mathsf{P}_{t,x}$. Using that $U^c(t, c_i(t)) = G(t, c_i(t))$ for all $t \in [t_*, T]$ [since $c_1$ and $c_2$ solve the system (4.7)–(4.8) as pointed out above] and that $U^c(T,x) = G(T,x)$ for all $x \geq 0$, we see that $U^c(t+\sigma_c, X_{t+\sigma_c}) = G(t+\sigma_c, X_{t+\sigma_c})$. Hence from (4.56) and (4.57) using the optional sampling theorem we find

$$\begin{aligned}U^c(t,x) &= \mathsf{E}_{t,x}(U^c(t+\sigma_c, X_{t+\sigma_c})) \\ &\quad - \mathsf{E}_{t,x}\bigg(\int_0^{\sigma_c} H(t+u, X_{t+u}) \\ &\qquad\qquad\qquad \times I(c_1(t+u) < X_{t+u} < c_2(t+u))\,du\bigg) \\ &= \mathsf{E}_{t,x}(G(t+\sigma_c, X_{t+\sigma_c})) \\ &\quad - \mathsf{E}_{t,x}\bigg(\int_0^{\sigma_c} H(t+u, X_{t+u})\,du\bigg) = G(t,x)\end{aligned}$$
(4.59)

since $X_{t+u} \in (c_1(t+u), c_2(t+u))$ for all $u \in [0, \sigma_c)$. This proves that $U^c(t,x) = G(t,x)$ for all $x \in [c_1(t), c_2(t)]$ with $t \in [t_*, T]$ as claimed.

10. We show that $U^c(t,x) \geq V(t,x)$ for all $(t,x) \in [0,T] \times [0,\infty)$. For this, consider the stopping time

(4.60) $\quad \tau_c = \inf\{0 \leq s \leq T-t \mid c_1(t+s) \leq X_{t+s} \leq c_2(t+s)\}$

under $\mathsf{P}_{t,x}$ with $(t,x) \in [0,T] \times [0,\infty)$ given and fixed. The same arguments as those following (4.58) above show that $U^c(t+\tau_c, X_{t+\tau_c}) = G(t+\tau_c, X_{t+\tau_c})$. Inserting $\tau_c$ instead of $s$ in (4.56) and using the optional sampling theorem, we get

(4.61) $\quad \begin{aligned}U^c(t,x) &= \mathsf{E}_{t,x}(U^c(t+\tau_c, X_{t+\tau_c})) \\ &= \mathsf{E}_{t,x}(G(t+\tau_c, X_{t+\tau_c})) \geq V(t,x)\end{aligned}$

proving the claim.

11. We show that $c_1 \leq b_1$ and $c_2 \geq b_2$ on $[t_*, T]$. For this, suppose that there exists $t \in [t_*, T)$ such that $c_2(t) < b_2(t)$ and examine first the case when $c_2(t) > b_1(t)$. Choose a point $x \in (b_1(t) \vee c_1(t), c_2(t)]$ and consider the stopping time

$$(4.62) \quad \sigma_b = \inf\{0 \leq s \leq T - t \mid X_{t+s} \leq b_1(t+s) \text{ or } X_{t+s} \geq b_2(t+s)\}$$

under $\mathsf{P}_{t,x}$. Inserting $\sigma_b$ in the place of $s$ in (4.51) and (4.56) and using the optional sampling theorem, we get

$$(4.63) \quad \begin{aligned} &\mathsf{E}_{t,x}(V(t+\sigma_b, X_{t+\sigma_b})) \\ &= V(t,x) + \mathsf{E}_{t,x}\left(\int_0^{\sigma_b} H(t+u, X_{t+u})\,du\right), \end{aligned}$$

$$(4.64) \quad \begin{aligned} &\mathsf{E}_{t,x}(U^c(t+\sigma_b, X_{t+\sigma_b})) \\ &= U^c(t,x) + \mathsf{E}_{t,x}\bigg(\int_0^{\sigma_b} H(t+u, X_{t+u}) \\ &\qquad\qquad\qquad \times I(c_1(t+u) < X_{t+u} < c_2(t+u))\,du\bigg). \end{aligned}$$

Since $U^c \geq V$ and $V(t,x) = U^c(t,x) = G(t,x)$ for $x \in [b_1(t) \vee c_1(t), b_2(t) \wedge c_2(t)]$ with $t \in [t_*, T]$, it follows from (4.63) and (4.64) that

$$(4.65) \quad \begin{aligned} &\mathsf{E}_{t,x}\left(\int_0^{\sigma_b} H(t+u, X_{t+u})I(X_{t+u} \leq c_1(t+u) \text{ or } X_{t+u} \geq c_2(t+u))\,du\right) \\ &\leq 0. \end{aligned}$$

Due to the fact that $H(t+u, X_{t+u}) > 0$ for $u \in [0, \sigma_b)$ we see by the continuity of $b_i$ and $c_i$ for $i = 1, 2$ that (4.65) is not possible. Thus under $c_2(t) < b_2(t)$ we cannot have $c_2(t) > b_1(t)$. If, however, $c_2(t) \leq b_1(t)$, then due to the facts that $b_1$ is decreasing with $b_1(T) = 0$ and $c_2(T) > 0$ there must exist $u \in (t, T)$ such that $c_2(u) \in (b_1(u), b_2(u))$. Applying then the preceding arguments at time $u$ instead of time $t$, we again arrive at a contradiction. Hence we can conclude that $c_2(t) \geq b_2(t)$ for all $t \in [t_*, T]$. In exactly the same way (or by symmetry) one can derive that $c_1(t) \leq b_1(t)$ for $t \in [t_*, T]$ completing the proof of the initial claim.

12. We show that $c_1$ must be equal to $b_1$ and $c_2$ must be equal to $b_2$. For this, let us assume that there exists $t \in [t_*, T)$ such that $c_1(t) < b_1(t)$ or $c_2(t) > b_2(t)$. Pick an arbitrary point $x$ from $(c_1(t), b_1(t))$ or $(b_2(t), c_2(t))$ and consider the stopping time $\tau_D$ from (3.1) under $\mathsf{P}_{t,x}$. Inserting $\tau_D$ instead of $s$ in (4.51) and (4.56), and using the optional sampling theorem, we get

$$(4.66) \quad \mathsf{E}_{t,x}(G(t+\tau_D, X_{t+\tau_D})) = V(t,x),$$

$$(4.67) \quad \begin{aligned} &\mathsf{E}_{t,x}(G(t+\tau_D, X_{t+\tau_D})) \\ &= U^c(t,x) + \mathsf{E}_{t,x}\left(\int_0^{\tau_D} H(t+u, X_{t+u})\right. \end{aligned}$$

$$\times I(c_1(t+u) < X_{t+u} < c_2(t+u))\,du\bigg),$$

where we also use that $V(t+\tau_D, X_{t+\tau_D}) = U^c(t+\tau_D, X_{t+\tau_D}) = G(t+\tau_D, X_{t+\tau_D})$ upon recalling that $c_1 \le b_1$ and $c_2 \ge b_2$, and $U^c = G$ either between $c_1$ and $c_2$ or at $T$. Since $U^c \ge V$ we see from (4.66) and (4.67) that

$$(4.68) \quad \mathsf{E}_{t,x}\bigg(\int_0^{\tau_D} H(t+u, X_{t+u}) I(c_1(t+u) < X_{t+u} < c_2(t+u))\,du\bigg) \le 0.$$

Due to the fact that $H(t+u, X_{t+u}) > 0$ for $X_{t+u} \in (c_1(t+u), c_2(t+u))$ we see from (4.68) by the continuity of $b_i$ and $c_i$ for $i = 1, 2$ that such a point $(t, x)$ cannot exist. Thus $c_i$ must be equal to $b_i$ for $i = 1, 2$ and the proof is complete. $\square$

REMARK 4.2. The following simple method can be used to solve the system (4.7)–(4.8) numerically. Better methods are needed to achieve higher precision around the singularity point $t = T$ and to increase the speed of calculation. These issues are worthy of further consideration.

Set $t_k = kh$ for $k = 0, 1, \ldots, n$ where $h = T/n$ and denote [recalling (4.3) and (4.4) above for more explicit expressions]:

$$(4.69) \quad \begin{aligned} I(t, b_i(t)) &= J(t, b_i(t)) - G(t, b_i(t)) \\ &= \mathsf{E}_{b_i(t)}(G(T, X_{T-t})) - G(t, b_i(t)), \end{aligned}$$

$$(4.70) \quad \begin{aligned} &K(t, b_i(t), t+u, b_1(t+u), b_2(t+u)) \\ &\quad = \mathsf{E}_{b_i(t)}(H(t+u, X_{t+u}) I(b_1(t+u) < X_{t+u} < b_2(t+u))), \end{aligned}$$

for $i = 1, 2$. Note that $K$ always depends on both $b_1$ and $b_2$.

The following discrete approximation of the integral equations (4.7) and (4.8) is then valid:

$$(4.71) \quad I(t_k, b_i(t_k)) = \sum_{j=k}^{n-1} K(t_k, b_i(t_k), t_{j+1}, b_1(t_{j+1}), b_2(t_{j+1}))h$$

for $k = 0, 1, \ldots, n-1$ where $i = 1, 2$. Setting $k = n-1$ with $b_1(t_n) = 0$ and $b_2(t_n) = 1/2\mu$, we can solve the system (4.71) for $i = 1, 2$ numerically and get numbers $b_1(t_{n-1})$ and $b_2(t_{n-1})$. Setting $k = n-2$ and using values $b_1(t_{n-1})$, $b_1(t_n)$, $b_2(t_{n-1})$, $b_2(t_n)$, we can solve (4.71) numerically and get numbers $b_1(t_{n-2})$ and $b_2(t_{n-2})$. Continuing the recursion we obtain $b_i(t_n)$, $b_i(t_{n-1})$ ,..., $b_i(t_1)$, $b_i(t_0)$ as an approximation of the optimal boundary $b_i$ at the points $T$, $T-h, \ldots, h, 0$ for $i = 1, 2$ (see Figures 1 and 2). Equation (4.11) can be treated analogously (see Figures 3 and 4).

PROGRAMME FOR ADVANCED MATHEMATICS OF FINANCE
SCHOOL OF COMPUTATIONAL AND APPLIED MATHEMATICS
UNIVERSITY OF THE WITWATERSRAND
PRIVATE BAG 3, WITS, 2050
SOUTH AFRICA
E-MAIL: jacques.du-toit@postgrad.manchester.ac.uk

SCHOOL OF MATHEMATICS
THE UNIVERSITY OF MANCHESTER
SACKVILLE STREET
MANCHESTER M60 1QD
UNITED KINGDOM
E-MAIL: goran@maths.man.ac.uk
URL: http://www.maths.man.ac.uk/goran